\documentclass[a4paper,12pt]{amsart}

\def\myfnt{\ifx\protect\@typeset@protect\expandafter\footnote\else\expandafter\@gobble\fi}
\makeatother

\usepackage{amsfonts}
 \usepackage{amssymb}
 \usepackage{amsmath,amsxtra,amsthm}

\usepackage[usenames]{color}

\usepackage[mathscr]{eucal}

 \usepackage{dsfont}

\usepackage{hyperref}

\usepackage{wrapfig}

\usepackage[all]{xy}

\usepackage{shuffle}
\usepackage{scalerel}[2016/12/29]

\newtheorem{theorem}{Theorem}
\newtheorem{deff}{Definition}

\newtheorem{example}{Example}
\newtheorem{lemma}{Lemma}

\newtheorem{cor}{Corollary}
\newtheorem{prop}{Proposition}
\newtheorem{rem}{Remark}

\renewcommand{\proof}{{\bf Proof.}~}

\newcommand{\mto}{\mapsto}
\newcommand{\bqa}{\begin{eqnarray}}
\newcommand\eqa {\end{eqnarray}}
\newcommand{\beq}{\begin{eqnarray}}
\newcommand{\beqn}{\begin{eqnarray}\nonumber}
\newcommand{\eeq}{\end{eqnarray}}
\newcommand{\be}{\begin{array}}
\newcommand{\ee}{\end{array}}

 \newcommand{\pt}{\partial}

   \newcommand\vf\varphi

 \newcommand{\Hom}{\mathrm{Hom}}

 \newcommand{\Id}{\mathrm{Id}}

 \newcommand{\uHom}{\underline{\mathrm{Hom}}}

 \newcommand{\cV}{{\mathcal V}}
 \newcommand{\cH}{{\mathcal H}}
 \newcommand{\cP}{{\mathcal{P}}}

 \newcommand{\cO}{{\mathcal{O}}}

 \newcommand{\cL}{\mathcal{L}}
 
 \newcommand{\cG}{{\mathcal{G}}}
 
 \newcommand{\cF}{{\mathcal{F}}}

\newcommand{\cA}{{\mathcal A}}

\newcommand{\cR}{{\mathcal R}}
\newcommand{\cS}{{\mathcal S}}

\newcommand{\oM}{{\overline{M}}}

 \newcommand{\ca}{\mathpzc{a}}
 
 \newcommand{\cs}{\mathpzc{s}}

 \newcommand{\cp}{\mathpzc{p}}

 \newcommand{\cz}{\mathpzc{z}}

 \newcommand{\C}{{\mathbb C}}
 
 \newcommand{\R}{{\mathbb R}}
 \newcommand{\Z}{{\mathbb Z}}

 \newcommand{\N}{{\mathbb N}}

\newcommand{\ua}{{\underline{a}}}
\newcommand{\ub}{{\underline{b}}}
\newcommand{\uU}{{\underline{U}}}
\newcommand{\uW}{{\underline{W}}}
\newcommand{\uphi}{{\underline{\phi}}}

  \newcommand{\noi}{{\vskip 2mm\noindent}}

   \def\a{\alpha}
   \def\b{\beta}

   \def\ga{\gamma}

 \def\bk{\mathds{k}}

\DeclareFontFamily{OT1}{pzc}{}
\DeclareFontShape{OT1}{pzc}{m}{it}{<-> s * [1.15] pzcmi7t}{}
\DeclareMathAlphabet{\mathpzc}{OT1}{pzc}{m}{it}

\newcommand{\cc}{\mathpzc{c}}

\newcommand{\Sym}{\mathrm{Sym}}
\newcommand{\sE}{\mathscr{E}}

\newcommand{\sV}{\mathscr{V}}

\newcommand{\sC}{\mathscr{C}}

\newcommand{\sI}{\mathscr{I}}
\newcommand{\sA}{\mathscr{A}}

\begin{document}

\bibliographystyle{amsplain}


\title[The category of  $\Z-$graded manifolds]{The category of $\Z-$graded manifolds: \\ what happens if you do not stay positive}
\author{Alexei Kotov}
\address{Alexei Kotov: Faculty of Science, University of Hradec Kralove, Rokitanskeho 62, Hradec Kralove
50003, Czech Republic}
\email{oleksii.kotov@uhk.cz}

\author{Vladimir Salnikov}
\address{Vladimir Salnikov: LaSIE  -- CNRS \&  La Rochelle University,
Av. Michel Cr\'epeau, 17042 La Rochelle Cedex 1, France}
\email{vladimir.salnikov@univ-lr.fr}

\begin{abstract} 
In this paper we discuss the categorical properties of $\mathbb Z$-graded manifolds. We start by describing the local model paying special attention to the differences in comparison to the $\N$-graded case. In particular we explain the origin of formality for the functional space and spell-out the structure of the power series. Then we make this construction intrinsic using filtrations. This sums up to proper definitions of objects and morphisms in the category. 
We also formulate the analogue of Batchelor's theorem for the global structure of $\Z$-graded manifolds. \\[-2em]
\end{abstract}

\maketitle


\renewcommand{\theequation}{\thesection.\arabic{equation}}

\vspace{-1em}
\section{Introduction / Motivation}\label{sec:introduction}

\noi This paper is a part of a series of works on the relation of graded geometry to theoretical physics and mechanics. For physics, It has been observed since a couple of decades that the language of differential graded manifolds and $Q$-bundles can be conveniently used in the context of studying symmetries and gauging of functionals: \cite{KS, SaSt, KSS, KS2}
 In mechanics, graded and generalized geometry can be used to spell-out the intrinsic structure of equations governing the evolution of a large class of systems. This can be then reflected in the choice or design of appropriate simulation tools: \cite{SH, RSHD}. 

\noi
On top of the mentioned ``applied'' problems we are necessarily interested in purely mathematical questions studying the graded analogues of classical geometric and analytic results. This particular text is in a sense a prequel to our recent paper \cite{DGLG} on the problem of integration of differential graded Lie algebras.  There, we have described the categories of graded Lie algebras, Lie groups and Harish-Chandra pairs, the latter one being an important step for finding the integrating objects; and then equipped all those with an appropriate differential structure. The main results of \cite{DGLG} were formulated for $\N$-graded, i.e. non-negatively $\Z$-graded Lie algebras. And writing it we realized that most of them hold true for the honest $\Z$-graded case, at the expense of properly defining the necessary objects in the framework. Since the subject is of proper interest, not necessarily related to integration, and some facts have non-trivial proofs, we decided to make a separate text out of it, hence this paper. 

\noi The presentation is organized as follows. We start by recalling the standard definition of the manifolds graded by a monoid, which we specify to the $\Z$-graded case. Then we describe the local model and in particular explain the unavoidable appearance of power series in non-zero degree variables -- this appearance is done ``step-by-step'' and results in the description of the structure of the functional space. Afterwords we revisit the construction in a more intrinsic algebraic form using filtrations of graded modules. We end up by the globalization result and formulate the $\Z$-graded analogue of the Batchelor's theorem.

\section{The local model} 
\label{section:gman}

\noi
Traditionally, to define a graded manifold ${M}$ one needs to provide the following data: $X$ -- the base manifold (the degree $0$ part, sometimes also called body), and $\mathcal{O}(M)$ -- the sheaf of graded functions,\footnote{A remark about notations: calligraphic letters will usually be related to graded functional spaces, while straight letters denote either smooth (degree zero) objects or ingredients of the graded ones.} which can be smooth, real or complex analytic, depending on the category we are dealing with.


\begin{deff}\label{def:general_grading} A grading is a mapping $\deg$ from local functions in $\mathcal{O}(M)$ to a commutative monoid $\Gamma$.  The domain of definition of $\deg$ corresponds to homogeneous functions. 
\end{deff}

\noindent Throughout the whole paper we assume that the commutation relations are governed by this grading, i.e. there exists a \emph{commutation factor} $\varepsilon \colon \Gamma^2 \to \R^{\times}$, such that   
for any couple of homogeneous graded objects $a$ and $b$, $a\cdot b = (-1)^{\varepsilon(deg(a), deg(b))} b \cdot a$.

\begin{deff}\label{deff:gradman_general}
A graded manifold $M$ is a topological space for which $\mathcal{O}(M)$ is locally modelled as functions on $(\uU, V)$, where $\uU$ is an open chart in a manifold\footnote{One can consider the base manifolds in smooth, real or complex analytic categories, unless the contrary is explicitly stated.}  $X$ and $V$ is a vector space graded in the above sense (we will also require that $V_0=\{0\}$). We call $U = (\uU, V)$ a graded chart of $M$. 
\end{deff}

\begin{rem}
We assume $M_0$ to be smooth, although most of what is presented holds for analytic manifolds or algebraic varieties. In Section \ref{sec:global} we will provide a couple of examples of what may fail in the complex analytic case.
\end{rem}
\begin{rem}
 For the degree of linear maps to be well-defined, the monoid $\Gamma$ needs to be cancellative, this reduces $\Gamma$ to combinations of $\N$, $\Z$, $\Z/p\Z$. 
\end{rem}

\noindent In this paper we mainly consider $\Gamma = \Z$, and thus $\varepsilon(a, b) = (-1)^{deg(a)deg(b)}$. 
\subsection*{Finite degree case -- minimal functional space}
\label{sec:series}

~ \\
In this section we will describe the construction of $\mathcal{O}(M)$ in the case of $\Z$-grading, and in particular explain why it is natural.

\noi 
Recall that in the classical case, when ${M}$ is $\N$-graded, 
homogeneous (w.r.t. the grading) functions in $\mathcal{O}(M)$ are polynomials in positively graded variables with coefficients  smoothly depending on degree $0$ variables\footnote{We will write ``degree'' for the value of the $\Z$-grading, in contrast to ``polynomial degree'' when it is appropriate and does not lead to a confusion.}, the whole $\mathcal{O}(M)$ contains finite sums of those. And the coordinate changes have the same form as homogeneous functions. We will use this case as a starting point for the general construction.

\noi
Let us now qualitatively understand what is the main difference of the $\N$-  and $\Z$-graded cases. In the $\N$-graded case one can in principal consider formal power series in  graded variables, but since all the degrees are positive the polynomial degree of homogeneous functions  has to be bounded. Now, we add negatively graded generators to the sheaf of functions, thus permitting non-trivial degree-zero combinations. It means that even making a polynomial coordinate change we potentially evaluate smooth functions of non-zero degree arguments, automatically ending up in formal power series in them that are obtained from the Taylor expansion. 

\noi
Note that generators of odd degree, unless their number is infinite, do not produce power series: since they are self-anticommuting,  their squares vanish, which again gives an upper bound on the polynomial degree. The even degree variables on the contrary have to be considered carefully. Let us make this statement more explicit. 
  Fix an open chart $\uU \subset \R^n$ of $X$, denote the coordinates there by $x$;  consider a (non-trivially) graded vector space ${V}$ and decompose it in the following way:
  \begin{equation} \label{sumV}
  {V} =\bigoplus_{i\in \Z\setminus 0} V_i^{d_i}= \ldots\oplus V^{d_{-l}}_{-l} \oplus V^{d_{-l+1}}_{-l+1}\oplus \dots \oplus V^{d_{-1}}_{-1} \oplus \{0\} \oplus V^{d_{1}}_{1} \oplus \dots \oplus V^{d_{k}}_{k}\oplus\ldots,
  \end{equation}
where  $V$'s are vector spaces, the subscripts of $V_{\bullet}^{\bullet}$ denote the degree of elements of the respective subspace, and  the superscripts -- the dimension of it. Add to it $V_0 \equiv \uU$ to obtain ${U} = (V_0 \oplus V) = (\uU, V)$.

\begin{deff} \label{def-finite-gr}
We say that the graded vector space $V$ (resp. graded manifold $M$) is
\begin{itemize}
\item  {of finite degree} if $|l| < \infty$, i.e.  
 the maximal/minimal degree of generating elements is bounded and the decomposition \ref{sumV}) indeed stops in both directions after a finite number of terms. 
 \item {of finite graded dimension} if $d_i = dim(V_i) < \infty, \forall i \in \mathbb Z$.
 \item {of finite dimension} if it is of finite graded dimension and of finite degree.
\end{itemize}

\end{deff}
For the moment and until the contrary is explicitly stated we assume the graded graded manifold to be of finite dimension.
\begin{rem}
Some subtleties will occur while discussing the space of maps between graded manifolds -- to obtain a closed category one will clearly have to relax this assumption.
\end{rem}

    

 \noi
Let us now describe the  space $ \mathcal{F}({U})$ of the most general graded functions on  ${U}$. 
 Consider two families of indices $i$'s and $j$'s  to distinguish between even and odd 
  degrees respectively, since only the parity of the element (not the degree) plays a role in commutation relations. For the moment we do not care about their ordering, just enumerate them. 
   Denote  for convenience
  \beqn 
  D_1 &=& d_{-2} + d_{-4} + \dots + d_2 + d_4 + \dots = \sum d_i, \\ \nonumber 
  D_2 &=& d_{-1} + d_{-3} + \dots + d_1 + d_3 + \dots = \sum d_j,
  \eeq
   respectively ``even'' and ``odd'' rank of ${M}$.  
   Decompose the set of variable into the subsets of 
   even and odd variables $(\zeta_1, \ldots, \zeta_{D_1})$ and $(\eta_1, \ldots, \eta_{D_2})$, respectively; the odd variables ($\eta$'s) are self-anticommuting, and thus square to zero,
   while the even ones ($\zeta$'s) are self-commuting and can be raised to arbitrary power. 
  Then a homogeneous $f$ of degree $k$ (all such functions constitute $\mathcal{F}({U})_k$) expands as a formal power series
   \begin{equation} \label{series}
    f = \sum\limits_{\begin{array}{c}
                     i_1, \dots, i_{D_1} \in \Z_{\ge 0} \\
                     j_1, \dots, j_{D_2} \in \Z_{2}
                    \end{array}} f_{i_1\dots i_{D_1}j_1\dots j_{D_2}}(x) \zeta_1^{i_1}\dots \zeta_{D_1}^{i_{D_1}} \eta_1^{j_1}\dots \eta_{D_2}^{j_{D_2}},
   \end{equation}  
 such that
 $$
 \sum_{\a=1}^{D_1} i_\a \deg \big( \zeta_\a\big)+\sum_{\a=1}^{D_2} j_\a \deg \big( \eta_\a\big)=k\,.
 $$   
   The whole functional space $\cF (U)$ consists of finite sums of homogeneous functions:
   $$
   \mathcal{F}({U})=\bigoplus_k \cF (U)_k\,.
   $$
   Morally $\cF (U)$ is
    `` $  \left( \mathcal{C}({U_0}) \right)^{ |\Z_{\ge 0}|^{D_1} \cdot 2^{D_2}}$ '', that is an infinite (but obviously countable!)  line of smooth functions that are ordered lexicographically by $j_1\dots j_{D_1}i_1\dots i_{D_2}$. We write $\mathcal{C}$ to mimic $C^\infty$, but as mentioned before the construction works for most of reasonable categories of base manifolds and corresponding classes of functions. 
    
\begin{theorem}[minimal class of functions]\label{thm:minimal}
$\cF(U)$ is the minimal class of functions, containing polynomials on variables of non-zero degree with coefficients being smooth functions of the degree zero variables, which is stable under degree preserving morphisms of the same type.
\end{theorem}
\noindent\proof It is the corollary of the consequent series of propositions \ref{prop:F1_generated}, \ref{prop:F1_stability} and \ref{prop:graded-series} below. $\blacksquare$

\noi
Knowing how these functions are constructed we can actually say much more about the structure of the series (\ref{series}). Let us start with a kind of ``negative result'', motivating the whole construction.
\begin{prop}
The space $\mathcal{P} := \{$ Polynomial functions of variables from ${V}$ with smooth coefficients of $x \}$ is not stable under degree zero coordinate changes of the same type.
\end{prop}
\noindent\begin{proof}
  A change of the type 
  \begin{equation} \label{coord-change}
  x' = x + \xi^{n_\xi}\psi^{n_\psi}, \text{ where } n_\psi \deg(\psi) = - n_\xi \deg(\xi)   
  \end{equation}
   is homogeneous of degree zero. Plugging $x'$ to any smooth function and spelling out its Taylor expansion around $x$, one obtains a power series in $\xi$ and $\psi$ that need not stop for degree reasons if the degrees of them are even.
If at least one even generator is present for the positively and negatively graded parts of ${V}$, such a combination can be found.   $\blacksquare$    
\end{proof}

\noi The last proposition means that one has to extend the space of functions.

\noi
Let us define 
$\mathcal{F}_1({U}) = \{$ finite combination of formal power series, where coefficients are polynomial in non-zero degree variables, and the series depend only on finite degree zero combinations of generators$\}$. It is obvious that $\cP (U)\subset \cF_1 (U)$.

\begin{prop}\label{prop:F1_generated}
The class of functions $\mathcal{P}$, together with algebraic operations and  polynomial degree-preserving morphisms generates $\mathcal{F}_1$.
\end{prop}

\noindent\proof
The inclusion of the generated space of functions into $\mathcal{F}_1$ is straightforward. The other way around, for a function in $\mathcal{F}_1$, it is sufficient to consider a power series of a degree zero combination as a (formal) Taylor expansion of some function. And since the combination is finite and the number of different non-zero degree monomials is finite as well, one repeats the procedure and traces back the starting polynomial.

\noi To be more precise, let $f\in\cF_1$ be a degree zero power series of a finite set of monomials $z_\a (\zeta, \eta)$, $\a=1, \ldots, r$, whose coefficients smoothly depend on $x$, i.e.
$$
f (x, \zeta, \eta)=F \left(x, z_\a (\zeta, \eta) \right), 
\hspace{2mm} \mathrm{where}\hspace{2mm} F \left(x, z_\a  \right)
=\sum_{i_1\ldots i_r} F_{i_1\ldots i_r}(x) z_1^{i_1}\ldots z_r^{i_r}\,.
$$
By the classical Borel theorem, there exists a smooth function $h(x, z)$, such that $F$ is the Taylor expansion of $h$ with respect to $z$ at $0$. Consider a polynomial $\Z-$graded morphism $\phi$, such that $\phi^* (x)=x, \phi^* (z_\a)=z_\a (\zeta, \eta)$; then $f=\phi^* (h)$.
$\blacksquare$.

\begin{prop}\label{prop:F1_stability}
   $\mathcal{F}_1$ is stable with respect to natural (algebraic) operations, as well as to the coordinate changes of the same type.
   \end{prop}
\noindent\begin{proof}
The difference to the previous statement is that now formal power series are also permitted in the coordinate changes. However, they are only of degree zero, thus survive algebraic operations, and again arrange to a finite number of combinations ``indexed'' by the degrees of monomials they come from.  
$\blacksquare$
\end{proof}

\noi
The statement above shows that $\mathcal{F}_1({M})$ is the minimal sheaf of functions on a graded manifold ${M}$,  necessary to perform all the natural operations. It is obviously a subsheaf of $\mathcal{F}({M})$. In fact, a stronger statement holds.

\begin{prop} \label{prop:graded-series}
$\mathcal{F}_1({U})=\mathcal{F}({U})$, that is all the formality in power series in $\mathcal{F}({U})$ can be reduced to combinations of variables with the total degree $0$.

\end{prop}
\noindent\begin{proof}
One considers the most general form (\ref{series}) of a function in $\mathcal{F}({U})$. 
For simplicity of the presentation let us ignore the odd variables: they do not produce infinite sums. 
The goal is to show that one can always rearrange the terms to bring all the infinite series to degree zero parts. Let us split all variables from $V$ into positively graded $\xi$'s, negatively graded $\psi$'s and zero degree $x'$, then
 $$  f = \sum\limits_{I, J} f_{IJ}(x) \xi^I \psi^J,
 $$
where $I$ and $J$ are multi-indices, such that $\xi^I = \xi_1^{I_1}\ldots \xi_k^{I_k}$ and $\psi^J = \psi_1^{J_1}\ldots \psi_l^{J_l}$. 
Fixing some degree $c$ of a homogeneous function $f$ amounts to the condition
\begin{equation} \label{dioph}
   I \cdot \deg(\xi) + J \cdot \deg(\psi) = c,   
\end{equation}
where $\deg(\xi)$ and $\deg(\psi)$ are vectors of degrees of all the $\xi$'s and $\psi$'s respectively, and 
``$\cdot$'' denotes the scalar product. 
Let $\deg (\xi_i)=a_i$ and $\deg (\psi_j)=-b_j$, where all $a_i$ and $b_j$ are positive integers.
We have thus reduced the problem of rearranging the terms to a classical Diophantine equation  $\ua \cdot I - \ub \cdot J = c$ with positive coefficient vectors $\ua$ and $\ub$,  for which we consider non-negative solutions.  This turns out to be a well-studied question from the theory of semigroups (\cite{clifford1, clifford2}) --  a commutative monoid can be generated by a finite basis. For completeness of the presentation we will explain this fact\footnote{We noticed that this fact has been revisited in the context of computer algebra (cf. for example \cite{clausen}) to produce efficient solutions of \ref{dioph}. It may be related to the discussion in \cite{SH-ca}.} here in details. 

\noi Let $a_1, \ldots,a_k $, $b_1, \ldots, b_l$ and $c$ be positive integers, $S(\underline{a}, \underline{b}, c)$
be the set of all non-negative integer solutions $(p_1, \ldots, p_k, q_1, \ldots, q_l )$ to the inhomogeneous linear Diophantine equation
\beqn
a_1 p_1 + \ldots + a_k p_k - b_1 q_1 - \ldots - b_l q_l =c 
\eeq
and $S(\underline{a}, \underline{b})$ be the set of all non-negative integer solutions to the corresponding equation
\beqn
a_1 p_1 + \ldots + a_k p_k - b_1 q_1 - \ldots - b_l q_l =0
\eeq
Denote by $M(\underline{a}, \underline{b}, c)$ and $M(\underline{a}, \underline{b})$ the subsets of $S(\underline{a}, \underline{b}, c)$
and $S(\underline{a}, \underline{b})\setminus \{0\}$, respectively, which are minimal with respect to the partial ordering
\beqn
(\underline{p}, \underline{q}) \le (\underline{p}', \underline{q}')  \Longleftrightarrow \left( \forall i \colon p_i\le p_i'   \right) \wedge
 \left( \forall j \colon q_j\le q_j'   \right)
\eeq

\begin{lemma}[Finiteness properties for linear Diophantine equations]\label{lem:finiteness}
{\mbox{}\vskip 2mm}
\begin{enumerate}
\item The sets  $M(\underline{a}, \underline{b}, c)$ and $M(\underline{a}, \underline{b})$ are finite.
\vskip 1mm
\item  $S (\underline{a}, \underline{b})$ is the set of all $\N-$linear combinations of elements in $M(\underline{a}, \underline{b})$.
\vskip 1mm
\item $S (\underline{a}, \underline{b}, c)  = M(\underline{a}, \underline{b}, c) + S (\underline{a}, \underline{b})$. 
\end{enumerate}
\end{lemma}

\noi
Hence the solutions of (\ref{dioph}) behave as we expect them to from linear algebra: the general solution decomposes into a (finite) part from inhomogeneous equation and a (potentially infinite, but finitely generated) part from the homogeneous one. Every solution of the inhomogeneous part will give a degree $deg$ monomial, while each infinite  family from the homogeneous part produces a degree zero power series. 
Here it is important to recall that we assumed the graded manifold to be of finite degree and each homogeneous component of $V$ to be finite dimensional.
Putting all this together, in our terms it means that the series above decomposes into a finite sum of degree $0$ formal series, which are multiplied by degree $c$ polynomials. 

\noi Moreover, Lemma \ref{lem:finiteness} states more: those degree $0$ formal power series are power series of a finite number of degree $0$ monomials, corresponding to the set of multi-indices $M(\ua, \ub)$.

\noi Indeed, let $c$ be the degree of $f$ (without loss of generality we assume that $c$ is a non-negative integer). 
By Lemma \ref{lem:finiteness} we can decompose $f$ into the finite sum of degree $c$
components 
\beq\label{eq:f-series} f(x,\xi,\psi)=\sum_{(I,J)\in M(\underline{a}, \underline{b}, c)} f_{IJ}(x, z (\xi, \psi)) \xi^I \psi^J\,,
\eeq 
where, with the notations from Proposition \ref{prop:graded-series}, $\xi$ and $\psi$ are negative and positive degree variables, $f_{IJ}(x, z))$ are smooth functions of zero degree variables $x$ and formal power series of zero degree variables $(z)=\{ z_{IJ}\}$, parametrized by $(I ,J)\in M(\underline{a}, \underline{b})$,
such that $z_{IJ} (\xi, \psi)=\xi^I\psi^J$. This is precisely the desired result.   $\blacksquare$

\end{proof}

\begin{rem} In \cite{DGLG} we have already discussed the properties of the functional space on graded manifolds (cf. appendix there). Comparing with the results of the current paper, in \cite{DGLG} we have considered the most general (i.e. largest possible) functional space not to worry about the operations (algebraic ones and compositions) to respect the space -- this leads naturally to consideration of formal series, as observed in \cite{Fairon:1512.02810, vysoky} as well. In this paper we do the opposite -- consider the minimal space (or rather the space looking ``small'') and prove that the operations do not lead out of it. \end{rem}

\noindent The following theorem is a corollary of Proposition \ref{prop:graded-series}.

\begin{theorem}[$\Z-$graded Borel's lemma] \label{cor:gr-borel}
Any homogeneous local formal power series is a Taylor expansion of some local smooth function of the same degree. \end{theorem}
\noindent\proof
 By the classical Borel's lemma, with the notations from Proposition \ref{prop:graded-series}, formula \eqref{eq:f-series}, we are able to choose local smooth functions $h_{IJ}(x, z)$, the Taylor expansion of which
with respect to $z$ gives us $f_{IJ}(x, z)$. Then the desired smooth degree $c$ function is
\beqn h (x, \xi,\psi)\colon =\sum_{(I,J)\in M(\underline{a}, \underline{b}, c)} h_{IJ}(x, z (\xi, \psi)) \xi^I \psi^J
\eeq 
$\blacksquare$


\section{The local model: algebraic description} 
\label{section:local_algebraic}

\noi In the previous Section \ref{section:gman} we substantiate the appearance of graded formal power series as appropriate functions on the local model of a $\Z-$graded finite-dimensional manifold.\footnote{Manifolds of this type are natural to be called semi-formal as local functions on them are formal power series of non-zero degree coordinates; later on (in Section \ref{sec:global}) we will see that, in the smooth category, such a manifold is always a formal neighborhood of its base embedded into some smooth supermanifold. } Now we are going to describe the local model algebraically in a regular way, by dropping of the assumption about finite number of non-zero degree coordinates.

\noi
\subsection{Graded modules.}\label{sec:graded_modules}

Let $R$ be a ring. 
\begin{deff}[Graded modules]
An $R-$module graded by a monoid $\Gamma$ is the direct sum of a collection of $R-$modules parameterized by $\Gamma$:
\beqn
\sE = \bigoplus_{i\in\Gamma} \sE_i\,.
\eeq
We say that $\sE$ is of finite degree if only a finite subset of modules $\sE_i$ are not equal to zero. Assume that every $\sE_i$ is a free module of finite rank $d_i$, then $\sE$ is said to be of finite graded rank $(d_i)_{i\in\Gamma}$. If $d=\sum_{i\in\Gamma}d_i<\infty$, then $\sE$ is called a graded $R-$module of finite (global) rank $d$. If $R$ is a field $\bk$, the word ``rank'' is to be replaced with ``dimension''.
\end{deff}

\noi For two graded $R-$modules $\sE$ and $\sE'$, $\Hom (\sE, \sE')$ is the set of morphisms in the category of graded $R-$modules, i.e. degree preserving morphisms over $R$
$$
\Hom \big(\sE, \sE'\big) =\prod_{i\in\Gamma} \Hom \big(E_i, E'_i\big)\, .
$$

\noi Assume for simplicity that $\Gamma$ is an abelian group. For any graded $R-$module $\sE$ one can produce another graded module $\sE [k]$ with degree shifted by $k\in\Gamma$, such that
$\sE[k]_i\colon =\sE_{i+k}$. 
For fiber-linear coordinates the shift goes in the opposite direction: 
$$
(\sE)^*_{-i}=(\sE_i)^* \, \Longleftrightarrow \,
(\sE)^*=\bigoplus_i (\sE_i)^*[i]
$$
thus
$$
(\sE [k])^*=\bigoplus_i (\sE[k]_i)^*[i]=\bigoplus_i (\sE_{k+i})^*[k+i][-k]=
\left(\bigoplus_j (\sE_{j})^*[j]\right)[-k] = (\sE)^*[-k]
$$
$(\sE)^*_{-i}=(\sE_i)^*$ comes from the requirement for $\sE^*\otimes \sE\to R$ being of degree 0.
 \\
In this formalism, a graded $R-$module $\sE$ is a collection of degree zero modules $E_i$ parameterized by $i\in \Gamma$, such that $\sE_i = E_i [-i]$.
When $k=1$ the operation $\sE\mto \sE[1]$, denoted by $s$, is called {\it the suspension}. It is clear that $\sE [k]=s^k \sE$.

\noi Along with the categorical $\Hom$, we also define $\uHom$, the set of ungraded morphisms of bounded degree: every such a morphism is a finite sum of morphisms of pure degree $p$, which change the degree of elements by some  $p\in\Gamma$:
$$
\uHom \big(\sE, \sE'\big) =\bigoplus_{p\in\Gamma} \uHom \big(\sE, \sE'\big)_p\, , \hspace{3mm}
\uHom \big(\sE, \sE'\big)_p= \left(\prod_{i\in\Gamma} \Hom \big(E_i, E'_{i+p}\big)\right)[-p]\,.
$$
As one can see from the definition, in contrast to $\Hom \big(\sE, \sE'\big)$, which is always concentrated in degree $0$, $\uHom \big(\sE, \sE'\big)$ is $\Gamma-$graded itself.

\noi  

\noi Henceforward $R$ will be a commutative unital ring. Under this assumption $\uHom$ will produce a graded $R-$module out of two graded modules.  

\noi Let $\sE$ be an $R-$module graded by $\Gamma$. Denote by $T^k \sE$ its $k-$th tensor power over $R$
\beqn
T^k \sE = \underbrace{\sE \otimes \ldots \otimes \sE}_k \,.
 \eeq
 viewed as an $R-$module, and by $T( \sE)$ the direct sum of all tensors: $T(\sE)=\bigoplus_{k\ge 0} T^k \sE $. We define the symmetric powers of $\sE$ over $R$ as follows\footnote{All tensor, symmetric, and exterior products and homomorphisms of modules over (graded) commutative rings 
are taken in the correspondent categories, unless contrary is assumed. Therefore the notation $\otimes$ in this context automatically means $\otimes_{R}$ and $\Sym(\sE)$ is the same as $\Sym_{R} (\sE)$.}
$$
 \Sym (\sE) = \mathrm{T}( \sE) / \left<  \cz_1\otimes \cz_2 - \varepsilon(\deg(\cz_1), \deg(\cz_2)) \cz_2 \otimes \cz_1  \, | \, \cz_1, \cz_2\in \sE\right>\,.
$$
Hereafter $\Sym(\sE)$ will
 be regarded as a free graded commutative $R-$algebra, generated by $\sE$ (depending on our purposes, we will also view it as an $R-$coalgebra).

\noi The dual module $\sE^*$ is defined as $\uHom (\sE, R)$, so that 
\beqn
\sE^*_{i} = \Hom \left(E_{-i}, R\right)[-i]\,.
\eeq
It is clear that if $\sE$ is a free graded $R-$module of finite graded rank, then so is $\sE^*$.

\begin{example}
Let $A$ be a commutative $\bk-$algebra and $\sV=\bigoplus_{i\in\Gamma} \sV_i$ be a graded vector space over a field $\mathds{k}$, corresponding to a collection $\big(V_i\big)_{i\in\Gamma}$ of vector spaces over $\bk$; as it was already mentioned in the beginning of the section, we do not anymore assume that $\sV$ is of finite degree, however, we claim that it is of finite graded dimension, i.e. all homogeneous components $\sV_i$ are finite-dimensional vector spaces. Consider the free graded $A-$module, generated by $\sV$ (by construction it has finite graded rank\footnote{Recall the definition \ref{def-finite-gr}.}): 
\beq\label{eq:module_out_of_vector_space}
\sE = \bigoplus_{i\in\Gamma} \sE_i \, , \hspace{5mm}
\sE_i \colon = A\otimes_{\bk} \sV_i\,.
\eeq 
Then $\sE^*_{i} = \Hom\left(E_{-i}, A\right)[-i]=
A\otimes_\bk \left(V_{-i}\right)^*[-i]$
and for any two free modules $\sE$ and $\sE'$, corresponding to graded vector spaces $\sV$ and $\sV'$, respectively, one has
$$
\uHom (\sE, \sE')=\bigoplus_{p\in \Gamma}A\otimes_{\bk}\left( \prod_{i\in\Gamma} \Hom \left(V_i, V_{i+p}\right)\right)\,.
$$
\end{example}

\noi


\subsection{Motivating example.}\label{sec:motivation_of_filtrations} Consider $\bk [x]$, the algebra of polynomials of one variable with coefficients in $\bk$. It is a non-negatively $\Z-$graded vector space: $\sV=\bigoplus_{n\ge 0} \bk^n [x]$, where $\bk^n [x]$ are polynomials of the polynomial power $n$. We can view it as graded algebra and coalgebra at the same time; in particular, the comultiplication is given by
$$
\bk [x]\ni P(x) \mto P(x+y)\in \bk [x,y]\simeq \bk [x]\otimes \bk[y]\,.
$$
The dual vector space in the category of $\Z-$graded vector spaces is 
$$
\sV^* = \bigoplus_{m\ge 0} \sV_{-m}^*\, ,
\hspace{3mm} \sV_{-m}^* \colon = \Big(\bk^m [x]\Big)^*\,.
$$
Thus any element of the dual graded space is a finite linear combination of the basic vectors $(u_{-m})_{m\in\Z}$, introduced such that
$$
u_{-m} \left(x^n\right)=\left\{
\be{cc}
1, & n=m \\
0, & n\ne m
\ee
\right.
$$
However, as soon as we ``forget'' about polynomial powers, the dual space will change. For instance, the above mentioned graded dual space is not stable under non-homogeneous change of variables. Indeed, let $\phi\colon x'\mto x+1$, then 
$$
\phi^* \big(u_{-m}\big) \left(\big(x'\big)^n\right)=\left\{
\be{cc}
\left(
\be{c}
n\\ m
\ee
\right), & n\ge m \\
0, & n < m
\ee
\right. 
$$
Therefore $\phi^* \big(u_{-m}\big)$ can be identified with an infinite formal sum
\beq\label{eq:formal_sum}
\phi^* \big(u_{-m}\big)=\sum_{n\ge m} \left(
\be{c}
n\\ m
\ee
\right) \big(u'\big)_{-n}\,,
\eeq
where $\big(u'\big)_{-n}$ is defined in the same way as $u_{-m}$.

\noi Since $\sV = \varinjlim F_p\sV$, where $F_p\sV=\bigoplus_{n<p} \sV_n$, the ``ungraded'' dual vector space has to be
\beq
\label{eq:projlim_polyn}
\varprojlim \left(\sV^*\big/ \mathrm{Ann}\big(F_p\sV \big)\right) \simeq
\varprojlim \left( \bigoplus_{m=0}^p \sV^*_{-m} \right)=
\prod_{m\ge 0} \sV^*_{-m}\,,
\eeq
where $\mathrm{Ann}()$ is the annihilator of the corresponding spaces. The difference between the direct product and the direct sum is that in the first case the number of non-zero homogeneous components is not necessarily finite. Such a (generally) infinite collection of elements of pure degree can be combined into a formal series (cf.~\eqref{eq:formal_sum}). Moreover, the comultiplication $\Delta$ on $\bk[x]$ inherits the graded dual space with the structure of a commutative algebra, isomorphic to the polynomial algebra of a dual variable $t$, such that $u_{-m}=t^m/ m!$ Now $\mathrm{Ann}\big(F_p\sV \big)$ becomes an ideal, so that the projective limit \eqref{eq:projlim_polyn} is a commutative algebra, naturally isomorphic to the algebra of formal power series $\bk [[x]]$.

\noi Assume now that there are several (or even infinite) number of variables $(x_i)$ of positive degrees, such that for any $k>0$ there is only a finite number of variables of degree $k$. Then $\bk[x_1, \ldots]$ is a non-negatively graded vector space. Moreover, it is of finite graded dimension as the homogeneous component $\big(\bk^n [x_1\ldots]\big)_k$ is zero for $n\ge k$ and is finite-dimensional for $n\le k$. Therefore the graded dual to $\bk^n [x_1\ldots]$ is also of finite graded dimension (notice that this property has nothing to do with another grading with respect to the polynomial powers); however, as soon as the degrees of variables are allowed to be negative, this is no longer the case (this feature was already mentioned in Subsection \ref{sec:series}). 

\noi In the next subsection  we will see how to come to an analogue of formal power series (see the construction above) in a purely algebraic way.


\noi
\subsection{Filtrations, completions and limits.}\label{sec:limits}
From now on by grading we will mean $\Z-$grading. Besides the homogeneous component of degree $0$ of any graded $R-$module is to be zero, i.e.
$$\sE=\bigoplus_{i\in\Z\backslash\{0\}}\sE_i\,.
$$

\noi


\noi
 Let us introduce the following increasing (decreasing) filtration\footnote{For convenience we have recollected some notations and results on limits and filtrations in the appendix \ref{app-filtr} } of $\Sym(\sE^*)$ ($\Sym (\sE)$, respectively):
by definition, for all $p >0$, $F^p\Sym(\sE^*)$ is the graded ideal of the symmetric algebra, generated by elements of degree $\le -p$, while $F_p \Sym (\sE)$ is the graded coideal of the symmetric coalgebra, cogenerated by
elements of degree $<p$, i.e. 
\beqn
F_p \Sym (\sE) \colon = \{a\in \Sym (\sE) \, |\, \Delta (a)\in \Sym (\sE)\otimes\Sym (\sE)_{<p}\}\,,
\eeq
where $\Delta$ is the comultiplication in the algebra $\Sym (\sE)$.

\vskip 2mm\noindent Let $\cR$ be the graded projective limit of $\Sym(\sE^*)/F^p \Sym (\sE^*)$ for $p\to \infty$.\footnote{This construction is compatible with the definition of the algebra of functions on a $\Z-$graded manifold given in \cite{Felder:2012kn}.} More precisely,
\beq\label{eq:proj_lim_R}
\cR \colon = \bigoplus_{i\in \Z} \varprojlim \Big(  \Sym(\sE^*)/F^p \Sym (\sE^*)\Big)_i\,,
\eeq
where
$$
\Big(  \Sym(\sE^*)/F^p \Sym (\sE^*)\Big)_i = 
\Sym(\sE^*)_i / F^p \Sym (\sE^*)_i\,.
$$

\noi It is worth mentioning here that both filtrations  are intrinsically defined (in terms of the corresponding graded algebraic structures on $\Sym (\sE')$, where $\sE'$ is equal to $\sE$ or $\sE^*$, i.e. we do not need to choose the space of generators $\sE'$; this property will play a crucial role in the definition of a graded manifold because an arbitrary morphism of algebras will not preserve the space of generators. 
However, once $\sE'$ is fixed, it is possible to make the construction more explicit by use of the $\N^2$ grading on $\Sym (\sE')$ for a graded module $\sE'$ (here $\sE'$ is $\sE$ or $\sE^*$): let $\deg_+$ and $\deg_-$ be the positive and negative degrees of the symmetric expressions counted independently, such that $\deg =\deg_+-\deg_-$. 

\noi In more details, suppose that $a\in \Sym (\sE')$ is decomposable in the following sense: $a=a_+ a_-$, where 
$a_\pm \in \Sym (\sE'_\pm)$ for
\beqn
\sE'_\pm &=& \bigoplus_{\pm i>0}\sE'_i\, , 
\eeq
then $\deg_+ a \colon= \deg a_+$ and $\deg_-a \colon= -\deg a_-$. 

\begin{lemma}[Filtrations in terms of generators]\label{lem:filtrations_in_generators}
For any $p>0$
\beqn
F_p \Sym (\sE) &=& \{ a\in \Sym (\sE) \, |\, \deg_+a<p \}
=\bigoplus_{q=0}^{p-1} \Big(\Sym \big(\sE_+\big)\Big)_{q}\Sym \big(\sE_-\big)\\ \nonumber
F^p \Sym (\sE^*) &=& \{ \cp \in \Sym (\sE^*) \, |\, \deg_-\cp \ge p \}=\bigoplus_{q=p}^{\infty} \Sym \big(\sE^*_+\big)\Big(\Sym \big(\sE^*_-\big)\Big)_{-q}\,,
\eeq
or 
\beqn
F_p Sym (\sE)_i &=& \bigoplus_{q=\max\{0,i\}}^{p-1} \Big(\Sym \big(\sE_+\big)\Big)_{q}
\Big(\Sym \big(\sE_-\big)\Big)_{i-q}\\ \nonumber
F^p Sym (\sE^*)_i &=& \bigoplus_{q= \max\{i,p\}}^\infty \Big(\Sym \big(\sE_+\big)\Big)_{q-i}
\Big(\Sym \big(\sE_-\big)\Big)_{-q}\,,
\eeq
\end{lemma}
\noindent\proof Straightforward computation.
$\blacksquare$

\begin{cor}[of Lemma \ref{lem:filtrations_in_generators}]\label{cor:filtration_in_generators}
For any $p\ge 0$, $ F_p \Sym (\sE)$
    has finite graded dimension.
\end{cor}
\noindent\proof
Indeed, 
$$
\big( F_p  Sym^k \sE\big)_i = \bigoplus_{l=0}^k \Big( F_p\left(\Sym^l\sE_+\Sym^{k-l}\sE_-\right)\Big)_i=
\bigoplus_{l=0}^k\left(
\bigoplus_{q=\max\{0,i\}}^{p-1} 
\big(\Sym^l \sE_+\big)_{q}\big(\Sym^{k-l} \sE_-\big)_{i-q}\right)\,.
$$
At each fixed degree $i$ it is non-zero only for a finite number of polynomial powers $k$: the claim is that $l$ and $k-l$ can not exceed $q$ and $q-i$, respectively, thus $k\le 2q-i< 2p-i$. 

\noi On the other hand, $\sE$ is of finite graded rank, i.e. the rank of each $\sE_j$ is finite. Therefore, for all admissible $l$ and $q$ (there is a finite number of such pairs for given $p$ and $i$), the modules $ \big( \Sym^l \sE_+\big)_q$ and $\big(\Sym^{k-l} \sE_-\big)_{i-q}$ have finite ranks, since
\beqn
\big( \Sym^l \sE_+\big)_i = \left( \Sym^l \left( \bigoplus_{j=1}^q \sE_j \right)\right)_q\, \hspace{2mm}
\mathrm{and} \hspace{3mm}
\big(\Sym^{k-l} \sE_-\big)_{i-q}=\left( \Sym^{k-l} \left( \bigoplus_{j=-1}^{i-q} \sE_{j} \right)\right)_{i-q}
\eeq
then so is $\Big( F_p\left(\Sym^l\sE_+\Sym^{k-l}\sE_-\right)\Big)_i$.
Combining all arguments together, we see that $ \left(F_p \Sym (\sE)\right)_i$ has finite rank for all $i<p$ (and is equal to zero otherwise). 
$\blacksquare$

\begin{lemma}
One has the canonical isomorphism of graded $R-$modules
\beq\label{isom:projlimit}
\cR\simeq \Big(\Sym (\sE)\Big)^*=
\uHom\left( \Sym (\sE),   R\right)
\eeq
\end{lemma}
\begin{cor}
Since $\Sym(\sE)$ is a graded $R-$coalgebra with the standard comultiplication
$$
\Delta\colon \Sym (\sE)\to \Sym (\sE\oplus\sE)=\Sym (\sE)\otimes \Sym (\sE)\,,
\hspace{2mm} \Delta (v)=v\otimes 1+1 \otimes v\, \hspace{2mm} \forall v\in\sE\,,
$$
provided \eqref{isom:projlimit} is proved,
we immediately deduce that $\cR$ is an $R-$algebra.
\end{cor}
\noindent\proof[Of Lemma \ref{isom:projlimit}] 
First we notice that $\Sym (\sE)=\varinjlim F_p (\sE) =\bigcup_{p\ge 0} F_p \Sym (\sE)$. Hence
\beq\label{eq:dual_is_projlim}
\Big(\Sym (\sE)\Big)^*\simeq \varprojlim
\Big(F_p\Sym (\sE)\Big)^*
\eeq

\noi The canonical pairing $\Sym (\sE^*)\otimes \Sym (\sE)\to R$ and the inclusions $F_p \Sym (\sE)\hookrightarrow \Sym (\sE)$ and $\Sym (\sE^*)\hookrightarrow \big(\Sym(\sE)\big)^*$
induces a morphism of graded $R-$modules
\beq\label{eq:epimorphism}
\Sym (\sE^*) \to \Big( F_p \Sym (\sE)\Big)^*
\eeq

\noi By Lemma \ref{lem:filtrations_in_generators}, the kernel of \eqref{eq:epimorphism} is $F^p \Sym (\sE^*)$.
By Corollary \ref{cor:filtration_in_generators},
    $ F_p \Sym (\sE)$ is of finite graded rank, thus 
    $\big(\Sym (\sE^*)\big)_i \to \big( F_p \Sym (\sE)\big)^*_i$
    is an epimorphism for all $i$. 
    
\noi    Therefore one has the following isomorphism of $R-$modules
    \beqn
    \big( F_p \Sym (\sE)\big)^*_i\simeq 
    \big(\Sym(\sE^*)\big)_i / \big(F^p \Sym (\sE^*)\big)_i \,.
    \eeq
    
\noi and thus and isomorphism of graded $R-$modules 
    \beq\label{eq:isomorphism_of_quotients}
    \Big(F_p\Sym (\sE)\Big)^* \simeq \Sym (\sE^*)/ F^p \Sym (\sE^*)\,.
    \eeq
This and \eqref{eq:dual_is_projlim} concludes the proof. $\blacksquare$

\begin{rem}\label{rem:projlim}
{\mbox{}\vskip 2mm}
\begin{enumerate}
\item Being defined as a projective limit (eq. (\ref{eq:proj_lim_R})), $\cR$ naturally inherits the filtration from $Sym(\sE^*)$
\beqn
\ldots\subset F^{p+1}\cR \subset F^p \cR \subset\ldots \subset F^1\cR \subset \cR \,,
\eeq
so that $F^p\cR$ annihilates $F_p Sym (\sE)$ by use of the isomorphism \eqref{isom:projlimit}. Therefore 
$\cR / F^p \cR$ is canonically isomorphic to $\Big(F_p\Sym (\sE)\Big)^*$ as well as to $\Sym (\sE^*)/ F^p \Sym (\sE^*)$.
\vskip 2mm
\item Given that $F_1\Sym (\sE)=\Sym (\sE_-)$, one has
$ \cR / F^1 \cR \simeq \Sym (\sE_-^*)$. In particular, the quotient algebra $\cR / F^1 \cR$ is non-negatively generated.
\vskip 2mm
\item The algebra $\cR$ can be canonically identified with the subalgebra of the algebra of formal power series, consisting of finite sums of homogeneous components.
\end{enumerate}
\end{rem}


\noi\subsection{Functor from $\Z-$graded to $\N^2-$graded algebras.}\label{sec:functor_to_bigraded} Let $\cA$ be a $\Z-$graded unital supercommutative algebra over $R$. There is an $\N^2-$graded algebra canonically associated to $\cA$. 
 First, consider the canonical decreasing filtration of $\cA$ by ideals $F^q\cA$ defined in the same way as above: $F^q\cA$ is generated by all elements of degree $\le -q$ for $q\ge 0$. In particular, $F^0\cA=\cA$. Obviously, $\big(F^q\cA \big)\big(F^{q'}\cA\big)\subset F^{q+q'}\cA$.
 
 \begin{rem}\label{rem:first_quotient}
    Let $\cA$ be a $\Z$-graded commutative algebra, then in the above notations $\cA_+ = \cA/F^1\cA$  is a non-negatively graded algebra. Moreover the quotient of $\cA_+$ by the ideal generated by all elements of degree higher than zero is the degree zero commutative algebra. This determines a canonical functor from the category of $\Z$-graded commutative algebras to the category of algebras commutative in the usual sense.

 \end{rem}

\noi Let $\cA^{gr}$ be the associated $\N-$graded algebra, 
\beqn
\cA^{gr}=\bigoplus_{q\ge 0} \cA^{gr}_{[q]}\,, \hspace{3mm}  \cA^{gr}_{[q]} \colon = F^{q}\cA / F^{q+1} \cA
\eeq
This algebra is bi-graded, such that 
\beqn 
\cA^{gr}_{p,q}\colon = 
\left( \cA^{gr}_{[q]}\right)_{p-q}\,,
\eeq
where the subscript on the r.h.s. comes from the original $\Z-$grading. Let $a\in F^q\cA$ be an element of pure degree, representing a non-zero element of the quotient space $F^{q}\cA / F^{q+1}\cA$, then $\deg a \ge -q$. From $\deg a=p-q$ we immediately obtain $p\ge 0$. Therefore 
\beq\label{eq:direct_sum_A}
\cA^{gr}=\bigoplus_{p,q\ge 0} \cA^{gr}_{p,q}
\eeq 
is an $\N^2-$graded algebra. Since the construction of the ideals, quotient spaces and the bigrading is canonical, we get a functor from the category of $\Z-$graded commutative algebras to the category of $\N^2-$graded commutative algebras over $R$. There is a right inverse functor which associates a $\Z-$graded commutative algebras to an $\N^2-$graded one: we take the same associative algebra and declare that $\deg a=p-q$ for any element $a$ of bi-degree $(p,q).$

\noi The algebra $\cA^{gr}$ is filtered, such that 
\beqn
F^q \cA^{gr}=\bigoplus_{j=q}^\infty \cA^{gr}_{-,j}
\eeq
This filtration is descended from the one of $\cA$, therefore the functor from $\Z-$graded to $\N^2-$graded algebras, defined above, is compatible with the filtrations.

\noi Let $\overline{\cA}=\varprojlim \cA^{gr}/F^q\cA^{gr}$; the submodule $\overline{\cA}_i$ can be viewed either as a direct product of bi-homogeneous $R-$submodules
\beq\label{eq:direct_prod_A}
\overline{\cA}_i =\prod_{p=\max\{i,0\}}^\infty \cA^{gr}_{p,p-i}=\prod_{q=0}^\infty \cA^{gr}_{q+i,q}
\eeq
or as an infinite power series, once we take the formal sum of the corresponding components. 
\noi
\begin{rem}\label{rem:bigraded_of_bigraded_and_completed}
If we apply the bi-grading functor $\cA\mto \cA^{gr}$ to $\overline{\cA}$, we obtain again $\cA^{gr}$, i.e. there is a canonical isomorphism
\beq\label{eq:bigraded_of_bigraded_and_completed}
\left(\overline{\cA}\right)^{gr}_{[q]} =
F^q \overline{\cA}/ F^{q+1} \overline{\cA} \simeq \cA^{gr}_{[q]} =
F^q \overline{\cA}/ F^{q+1} \overline{\cA}
\eeq
\end{rem}
\noi
\begin{rem}\label{rem:generalized_bigrading}
This algebra is bi-graded in the following (more general) sense: the two grading operators, being extended by the Leibniz rule to a pair of commuting Euler derivations $\varepsilon_{\pm}$, are acting on $\overline{\cA}$, such that
$$
\left(\varepsilon_{+}\right)_{|\cA^{gr}_{p,q} }=p\Id\,, \hspace{3mm} \left(\varepsilon_{-}\right)_{|\cA^{gr}_{p,q} }=q\Id\,.
$$
The total grading, which is descended from the original grading of $\cA$, is determined by the derivation $\varepsilon=\varepsilon_+-\varepsilon_-$. In contrast to
the direct sum \eqref{eq:direct_sum_A}, which consists of finite combinations of bi-homogeneous elements (or elements of pure bi-degree, i.e. joint eigenvectors of the Euler derivations $\varepsilon_\pm$), elements of the direct product \eqref{eq:direct_prod_A} are generally infinite sums of bi-homogeneous components. There exists a subalgebra of $\overline{\cA}$, the minimal subalgebra which contains all elements of pure bi-degree. 
\end{rem}

\begin{example}\label{ex:bigraded}
Consider the algebra $\cR$ as in Subsection \ref{sec:limits}. By Lemma \ref{lem:filtrations_in_generators},
$$
\cR^{gr}_{p,q} \simeq \Sym \big(\sE^*_+\big)_p  \Sym \big(\sE^*_-\big)_{-q}=
\bigoplus_{k=0}^{p}\bigoplus_{l=0}^{q}
\Big(\Sym^k \sE^*_+\Big)_p  \Big(\Sym^l \sE^*_-\Big)_{-q}
\,
$$
and
\beq\label{eq:direct_sum_R}
\cR^{gr}_i=\bigoplus_{p-q=i} \cR^{gr}_{p,q}=
\bigoplus_{p=\max\{i,0\}}^\infty \cR^{gr}_{p,p-i}
\,.
\eeq

\vskip 2mm\noindent The algebra $\cR^{gr}$ is filtered, such that for all $q\ge 0$
\beqn
F^q\cR^{gr}\colon = \bigoplus_{r=q}^\infty \cR^{gr}_{-,r}
\simeq \bigoplus_{r=q}^\infty \Sym \big(\sE^*_+\big)  \Sym \big(\sE^*_-\big)_{-r}\,.
\eeq
Let $\overline{\cR}=\varprojlim \cR^{gr}/F^q\cR^{gr}$ be constructed as above. 
Once we choose the module of homogeneous generators, $\overline{\cR}$ becomes isomorphic to $\cR$; furthermore, this induces an isomorphism of the ``finite bi-degree subalgebra'' $\overline{\cR}_{min}$ and $\Sym (\sE^*)$. 
\end{example}

\subsection{Morphisms of graded domains} 
\noi
In this subsection we describe the morphisms of graded domains. i.e. the mappings between them ``respecting'' the graded structure. We do it in a very explicit way using the local model, with the main goal to stress the difference between the $\N$- and $\Z$-graded cases. Those local descriptions can be glued together by a standard sheaf-theoretic argument, using the algebraic description we have defined above, that results then in the Definition \ref{def:gr-morph} of morphisms of graded manifolds.

\noi Recall that a graded coordinate chart $U$ is a couple $(\uU, V)$ consisting of a coordinate chart $\uU$ of the base manifold $X$, and a graded vector space $V = \bigoplus_{i\in \Z \setminus \{0\}} V_i$, such that $\cO_M(\uU)$ is locally modelled on $U$ in the sense of Definition \ref{deff:gradman_general}. We will denote the coordinates on $\uU$ by $x^i$, and on $V$ by $\theta^j$, for the degree of $\theta$'s which can now be  positive or negative. 
 
\noi Applying the procedure from the above Remark \ref{rem:first_quotient} to the algebra of graded functions on $U$ one canonically obtains the algebra of smooth functions on $\uU$. In particular any morphism of such graded algebras produces a morphism of corresponding algebras of smooth functions.

 \noi
\begin{deff} \label{def:gr-morph-loc}
A morphism of graded coordinate charts $\phi \colon U \to U'$ is a couple $\phi = ({\uphi}, \phi^\sharp)$, where $\phi^\sharp$ is a morphism of graded commutative algebras from the algebra of functions on $U'$ to the algebra of functions on $U$, i.e. a morphism of associative unital  algebras preserving the grading; and ${\uphi}$ is a smooth mapping between the underlying spaces such that $\uphi^*$ is obtained from $\phi^\sharp$ as above.
\end{deff}

\noi For the $\Z$-graded case, the definition mimics the $\N$-graded one, up to one important subtlety to be spelled out: the non-trivial relation of the mapping $\uphi \colon \uU \to \uU'$ and the degree zero restriction of the mapping  from $U$ to $U'$. The main difference, like for the functional spaces, is the possibility to construct a non-trivial degree zero combination out of non-trivially graded generators.
Let $U$ and $U'$ be two graded coordinate charts of $M$ and $M'$, with coordinates $(x^i, \theta^j)$ 
$(y^i, \tau^j)$, respectively.
The most general degree preserving mapping between $U$ and $U'$ reads: 
\begin{eqnarray}
\tau^i_{-k} &=& \sum_I f^i_{-k,I}(x) \theta^I, \quad I \cdot deg(\theta) = -k \label{neg} \\ 
y^i &=& f_0^i(x) + \sum_I f^i_{0,I}(x) \theta^I, \quad I  \cdot deg(\theta) = 0 \label{zer} \\
\tau_{k} &=& \sum_I f_{k,I}(x) \theta^I, \quad  I \cdot deg(\theta) = k, \label{pos} \\
&& k \le deg(M'),     \nonumber 
\end{eqnarray}
where in each line the summation is done over a multiindex $I$, and $deg(\theta)$ is the multiindex of degrees of the line of $\theta$'s, $deg(M')$ is the degree of the graded manifold $M'$, here supposed to be finite in the sense of the Definition \ref{def-finite-gr}.  

\begin{rem}
Regression to the $\N$-graded case. Clearly, if there are no negative degree generators on $U$, the first line (\ref{neg}) is absent, and the second line (\ref{zer}) stops after the first term. \end{rem}

\begin{rem}
It is important to note that all the sums in (\ref{neg} - \ref{pos}) are now formal power series in $\theta$'s, since the Diophantine condition $I \cdot deg(\theta) = d$ potentially admits infinitely many solutions. With the same argument of Lemma \ref{lem:finiteness}, one can reduce all the formality to degree zero. That means that for finite dimensional graded manifolds one can rewrite the infinite sums in \ref{neg} and \ref{pos} as finite sums of non-zero degree monomials with coefficients being power series of total degree zero.
\end{rem}

\begin{rem}
 The non-degeneracy of the mapping defined by (\ref{neg} - \ref{pos}) amounts to the non-degeneracy of the smooth maps $f^i_0$ in (\ref{zer}) and the non-degeneracy of the terms corresponding to degree $-k$ and $k$ generators in (\ref{neg}) and (\ref{pos}), which are linear.     
\end{rem}

\begin{prop}
A coordinate description (\ref{neg} - \ref{pos}) is enough to fully determine a morphism of graded domains.  
\end{prop}

\noi
\begin{proof}
 In the non-graded case this a classical result about the coordinate description of smooth mappings between smooth manifolds. For the $\N$-graded case, the proof is done by induction over the maximal degree of generators, starting from degree zero. Each induction step consists of adding polynomial functions in the appropriate degree, which is clear since their coefficients remain of degree zero, and by a general statement, morphisms of polynomial functions over any ring are defined by the images of the generators. 
 
 \noi
 For the honest $\Z$-graded case this induction is done starting from all the non-negative (or equivalently non-positive) generators and adding step-by-step the remaining negative (respectively, positive) ones. When adding an odd generator, the step is a trivial statement about polynomial algebra; for an even one we use the analogue of the Theorem \ref{cor:gr-borel}. Alternatively, one can use the approach of filtrations similar to the proof of the $\Z$-graded generalization of Batchelor's theorem below (Theorem \ref{gr-batchelor}).   
  $\blacksquare$
\end{proof}




\section{The global theory}\label{sec:global}
\noi In Section \ref{section:local_algebraic}, given a graded $R-$module $\sE$, we constructed a canonical graded $R-$algebra $\cR$ out of $\Sym (\sE^*)$ by use of the algebraic structure and $\Z-$grading on $\Sym (\sE^*)$. Whenever the module of generators $\sE$ is fixed, $\cR$ becomes isomorphic to the graded algebra of formal power series on $\sE^*$ with coefficients in $R$; in this case, in addition, $\cR$ inherits an $\N^2-$grading. In subsection \ref{sec:functor_to_bigraded}, we associated to any $\Z-$graded algebra an $\N^2-$graded algebra in a natural (functorial) way. Now we would like to apply this algebraic technique to the regular description of $\Z-$graded manifolds.

\noi Hereinafter a ring will be the algebra of local (smooth, analytic or algebraic) functions on a topological space with the corresponding structure. 


\noi \subsection{$\Z-$graded manifolds.}\label{sec:manifold_gluing} It was already mentioned in Section \ref{section:gman}, Definition \ref{deff:gradman_general} that a $\Z-$graded manifold is a ringed space locally modelled by functions 
on a graded coordinate chart $U = (\uU, V)$, where $\uU$ is an open coordinate chart and $V$ is a $\Z-$graded vector space, $V=\bigoplus_{i\in \Z\setminus \{0\}}V_i$. The precise way of local modelling, i.e. the definition of the algebra of graded functions on $U'= (\uU',V)$, where $\uU'$ is an open subset of $\uU$, is determined by formula \eqref{eq:proj_lim_R} in Subsection \ref{sec:limits}: here $R=\cO (\uU')$ and $\sE$ is the free $R-$module produced out of $\cV$ by use of \eqref{eq:module_out_of_vector_space}. Local graded functions, denoted by $\cF=\cF_U$, constitute a free $\cO_U-$module, thus it must be a sheaf over $\uU$. As a consequence, $\cF$ is a sheaf of $\cO_U-$algebras over $\uU$; this makes $\uU$ into a ringed space.  

\noi A global $\Z-$graded manifold is built from local ones in the standard way. Let $\{\uU_\alpha\}_{\alpha\in I}$ be an open cover of a topological (smooth, analytic or algebraic) space $X$ by (smooth, analytic or algebraic) coordinate charts and let $\psi_{\a\b}\colon \cF_{\uU_\b}\simeq\cF_{\uU_\a}$ be isomorphisms of the corresponding sheaves of graded functions over $\uU_{\a\b}=\uU_\a\cap \uU_\b$ (i.e. $\psi_{\a\b}$ are isomorphisms of sheaves, which respect the algebraic structure), satisfying cocyclic conditions over triple overlaps: $\psi_{\a\ga}=\psi_{\a\b}\psi_{\b\ga}$ for all $(\a,\b,\ga)$ such that $\uU_{\a\b\ga}=\uU_\a\cap \uU_\b\cap \uU_\ga \ne\emptyset$. 

\begin{deff}\label{def:sheaf_of_graded_functions}
For any open $\uU'\subset X$, let
\beqn 
\cF (\uU')=\{\big(f_\a\big)_{\a\in I}\in\prod_{\a\in I}\cF_{\uU_\a'}\,|\, (f_\a)_{| \uU_{\a\b}'} = \psi_{\a\b}(f_\b)_{| \uU_{\a\b}'}\, for\, all\,\, \uU_{\a\b}'\ne\emptyset\}\,,
\eeq
where $\uU_\a'=\uU'\cap \uU_\a$ and $\uU_{\a\b}'=\uU_\a'\cap \uU_\b'=\uU'\cap \uU_{\a\b}$.
\end{deff}

\noi By the standard argument from the sheaf theory, $\cF$ is a sheaf of $\Z-$graded algebras over $X$, i.e. it is a presheaf of algebras, such that for any local set $\uW\subset X$ and any open covering $(\uW_\a)_{\a\in I}$ of $\uW$, a collection of local sections $f_\a\in\cF(\uW_\a)$, $\a\in I$, compatible over double overlaps, uniquely determines a section $f\in\cF(\uW)$, such that $f_\a=f_{|\uW_\a}$. In other words
\beqn
\cF (\uW)\to \prod_{\a\in I} \cF(\uW_\a) \rightrightarrows \prod_{(\a,\b)\in I^2} \cF(\uW_{\a\b})
\eeq
is an equalizer.

\noi
It is now natural to give the global analogue of the definition \ref{def:gr-morph-loc}.
\begin{deff}  \label{def:gr-morph}
 A morphism of graded manifolds is a morphism of corresponding ringed spaces.
 If it is invertible it is called an isomorphism. 
\end{deff}
\begin{rem}
 By definition, the inverse to an isomorphism is obviously and isomorphism.
\end{rem}


\subsection{The Batchelor's theorem for smooth $\Z-$graded manifolds}

\vskip 2mm

\noindent A $\Z-$graded manifold $M$ is a topological (smooth, analytical) space $X$ together with a sheaf of $\Z-$graded functions, constructed as in subsection \ref{sec:manifold_gluing} for some open covering $(\uU_\a)_{\a\in I}$; from now on we will call it the structure sheaf of $M$ and denote by $\cO_M$. 

\begin{example}[$\Z-$graded vector bundles]\label{ex:graded_vector_bundle} 
Let $\cV=\bigoplus_{i\in\Z\setminus \{0\}}\cV_i$ be a $\Z-$graded vector bundle over $X$. Then the total space of $\cV$ admits a canonical structure of a $\Z-$graded manifold, such that the structure sheaf is determined by the procedure described in Subsection \ref{sec:limits}:
\beqn
\cO_{\cV}=\varprojlim \big(\Sym (\cV)/F^q \Sym (\cV)\big)\,,
\eeq
where $\Sym(\cV)$ is the sheaf of local sections of the bundle of symmetric powers of $\cV$ and $F^q$ is the filtration defined as in Subsection \ref{sec:limits}. Taking into account that the bundle of generators is fixed, there is a canonical isomorphism of sheaves of $\Z-$graded algebras
$\cO_{\cV}=\overline{\cO_{\cV}}$, so that
$\left(\overline{\cO_{\cV}}\right)_{min}=\Sym(\cV)$ (see Example \ref{ex:bigraded} for the explanation).
\end{example}

\noindent Given that $\cO_M$ is a sheaf of $\Z-$graded algebras, we endow it with the canonical decreasing filtration $F^p$, which was introduced in the Subsection \ref{sec:functor_to_bigraded}:
\beqn
\ldots\subset F^{p+1}\cO_M \subset F^p \cO_M \subset\ldots \subset F^1\cO_M \subset \cO_M \,,
\eeq
Notice that $\cO_M (\uU_\a)$ is isomorphic to a freely generated $\Z-$graded $\cO_{X} (\uU_\a)-$algebra  (in the sense of Subsection \ref{sec:limits}) for any coordinate chart $\uU_\a$ from the open covering. By use of Remark
\ref{rem:first_quotient}  
we conclude that $\cO_M/F^1\cO_M$ is a sheaf of non-negatively graded algebras. By Remark \ref{rem:projlim}, Statement 2, it is freely generated as an $\cO_X (\uU_\a)-$algebra over open coordinate charts $\uU_\a$, thus $F^1\cO_M$ can be regarded as the ideal of a (canonically defined) non-negatively graded submanifold $M_+$, such that in any local coordinates it is determined by the equation $\{z^i=0\}$ for all coordinates of degree $<0$. Likewise, there exists a canonical non-positively graded submanifold $M_-$ of $M$, the ideal of which is generated by all local functions of positive degrees. It is obvious, again in the spirit of the Remark \ref{rem:first_quotient}, that the zero-degree part (the body) of $M_\pm$ is $X=M_0$; moreover, $X=M_-\cap M_+$. 

\vskip 2mm\noindent Similar to the case of super manifolds, the base (the body) of a $\Z-$graded manifold is embedded into it, however, in contrast to $\N-$graded manifolds, in general there is no projection onto the base.  Recall that, whenever the grading is non-positive or non-negative, the projection is given by the inclusion of the structure sheaf of the base as the zero degree part. As soon as there are generators of mixed degrees, this procedure stops working. Nevertheless, by use of the technique described in Subsection \ref{sec:functor_to_bigraded}, we may construct a canonical graded manifold $\oM$, associated to $M$, the structure sheaf of which is $\cO_{\oM}\colon =\overline{\cO_M}$. The manifold $\oM$ is $\N^2-$graded (in the sense of Remark \ref{rem:generalized_bigrading}), thanks to the construction of its structure sheaf.

\begin{lemma}\label{lem:bar_iso}
For any smooth $\Z-$graded manifold $M$ there exists an isomorphism $M\simeq\oM$ in the category of $\Z-$graded manifolds.
\end{lemma}
\begin{rem}
 In this case the definition of an isomorphism rather straightforward. It can be spelled-out explicitly using the local model of section \ref{deff:gradman_general}, keeping in mind the non-degeneracy condition which starts at the degree $0$ level and then is extended by induction to positive and negative degrees. The same construction can be done in a more intrinsic way in the language of filtration of section \ref{section:local_algebraic} or as a regression from a more general construction of morphisms in the category of $\Z$-graded manifolds of the Definition \ref{def:gr-morph}.
\end{rem}
\noindent\proof We will construct such an isomorphism by induction, i.e. for any $p>0$ we will show that there exists an isomorphism of sheaves 
$\phi^p\colon \cO_M/F^p\cO_M\simeq \cO_\oM/ F^p \cO_\oM$, such that $\phi^{q+1}\mathrm{mod} \,F^q\cO_M=\phi^q$. Then the collection of $(\phi^p)$ will determine the desired isomorphism between the structure sheaves.

\vskip 2mm\noindent The first term $\phi^1$ exists by construction of $\cO_\oM$ (and it is canonical). Assume that we already have $\phi^p$. Let $(\uU_\a)_{\a\in I}$
be an open covering of $X$ by coordinate charts, which determines the structure of a $\Z-$graded manifold (as in Subsection \ref{sec:manifold_gluing}). Now one can extend $\phi^p$ to an isomorphism 
\beqn
\varphi^{p+1}_\a \colon \left(\cO_M/F^{p+1}\cO_M\right)(\uU_\a)\simeq \left(\cO_\oM/F^{p+1}\cO_\oM\right) (\uU_\a)
\eeq 
over each open chart $\uU_\a$ (given by the choice of local homogeneous coordinates). 

\vskip 2mm\noindent
Since $\varphi^{p+1}_\a $ and $\varphi^{p+1}_\b$ coincide modulo $F^p\cO_M$ over $\uU_{\a\b}$, we conclude that their difference 
\beqn
\psi_{\a\b}^{p+1}\colon = \varphi^{p+1}_\a -\varphi^{p+1}_\b\,,
\eeq 
being composed with the quotient map $\cO_\oM/F^{p+1}\cO_\oM\to \cO_\oM/F^{p}\cO_\oM$, gives us zero. Hence 
$\psi_{\a\b}^{p+1}$ takes values in the quotient $\cL^p (\uU_{\a\b})$, where
$$
\cL^p\colon=\left(\cO_\oM\right)^{gr}_{[p]} =\left(F^p \cO_\oM/F^{p+1}\cO_\oM\right)\subset \cO_\oM/F^{p+1}\cO_\oM\,.
$$
Consider the composition of morphisms
$$
\cO_M (\uU_\a) \to \left(\cO_M/F^{p+1}\cO_M\right)(\uU_\a) \xrightarrow{\varphi^{p+1}_{\a}} \left(\cO_\oM/F^{p+1}\cO_\oM\right) (\uU_{\a})
$$
which we will denote by the same letter. Using the multiplicative property of the filtration $\big( F^q \cO_\oM\big)\big( F^{q'} \cO_\oM\big)\subset F^{q+q'} \cO_\oM$, we conclude that
\beq 
im\big(\psi_{\a\b}^{p+1}\big) \big(\cL^p (\uU_{\a\b})\big) \subset \cL^p (\uU_{\a\b}) \cL^p (\uU_{\a\b})
\subset \cL^{2p} (\uU_{\a\b})
\eeq
For any $p>1$ we have $2p>p+1$, therefore $im (\psi_{\a\b}^{p+1}) \cL^p (\uU_{\a\b})$ is zero modulo $F^{p+1}\cO_\oM$. 
This implies that for any local $f$ in $\cO_M$ one has
\beqn
\big(\varphi^{p+1}_\a (f)-\varphi^{p+1}_\b (f)\big) \cL^p (\uU_{\a\b})=0
\eeq
and thus 
the action of $\cO_M$ on $\cL^p$ via $\varphi^{p+1}_\a$ and $\varphi^{p+1}_\b$ coincide over $\uU_{\a\b}$, which makes $\cL^{p}$ into a globally defined $\cO_M-$module.  Using this fact and the morphism property of $\varphi^{p+1}_\a$ for all $\a$, we conclude that $\psi_{\a\b}^{p+1}$ is a derivation of $\cO_M (\uU_{\a\b})$ with values in $\cL^p$; moreover, it is a \v{C}ech 1-cocycle. 

\noi Notice that $\cL^p$ is also a module over $\cO_{\oM_+}=\cO_{\oM}/ F^1\cO_{\oM}$ since $\big(F^1\cO_{\oM}\big) \big(F^p \cO_{\oM}\big)\subset F^{p+1} \cO_{\oM}$, therefore $F^1\cO_{\oM}$ annihilates $\cL^p$. Taking into account that $\cO_X$ is canonically embedded into the non-negatively graded sheaf of algebras $\cO_{\oM_+}$ as the subspace of degree $0$ elements, we conclude that $\cL^p$ is a well-defined $\cO_X-$module and so is the sheaf of derivations of $\cO_M$ with values in $\cL^p$. This allows to apply the partition of unity argument in the smooth category. Now we see that  $\psi_{\a\b}^{p+1}$ is necessarily a coboundary: there exists a collection of derivations $\psi^{p+1}_\a$ of the same type over all $\uU_\a$, satisfying \beqn\psi_{\a\b}^{p+1} =\psi^{p+1}_\b -\psi^{p+1}_\a\,.\eeq 

\vskip 2mm\noindent
Now we put
$\phi^{p+1}_\a \colon = \varphi^{p+1}_\a +\psi^{p+1}_\a$; it is again a morphism and, furthermore, two such morphisms coincide over $\uU_{\a\b}$, therefore we obtain a globally defined required morphism of sheaves.
$\blacksquare$


\vskip 2mm\noindent The Batchelor's theorem for smooth $\N-$graded manifolds asserts the existence of a canonical non-negatively graded vector bundle 
$\cV_+\to M_0$ together with a non-canonical isomorphism of $\N-$graded manifolds $M_\mp$ and the total space of $\cV_+$. The following theorem extends this result to the smooth $\Z-$graded case.

\begin{theorem}[$\Z$-graded Batchelor's theorem] \label{gr-batchelor}
There exists a non-canonical isomorphism of $\Z-$graded smooth manifolds between $M$ and the total space of $\cV=\cV_-\oplus \cV_+$.
\end{theorem}
\noindent\proof We apply Lemma \ref{lem:bar_iso} and the corresponding theorems for $\N-$graded manifolds.
Indeed, as soon as we fix an isomorphism $M\simeq\oM$ there are two canonical projections $M\to M_+$ and $M\to M_-$, given by the inclusions of the structure sheaves
\beqn
\cO_{M_\pm}\hookrightarrow \cO_{\oM}\simeq \cO_M \,.
\eeq
Now one has $M=M_+\times_{M_0} M_-$, where the fibered product of graded manifolds is defined algebraically in terms of the corresponding sheaves of functions: $\cO_M=\cO_{M_+}\otimes_{\cO_{M_0}}\cO_{M_-}$.
$\blacksquare$

\begin{example}[non-splitted $\N-$manifold]\label{ex:non-splittedN} 
Assume that $\cV=\cV_{-2}\oplus \cV_{-4}$ is a negatively graded vector bundle over $M_0$, where $\cV_{-2}$ and $\cV_{-4}$ are homogeneous summands of degree $-2$ and $-4$, respectively. As usual, by convention $\cV_{-2}=V[-2]$ and $\cV_{-4}=V[-4]$ for some degree $0$ vector bundles $V_2$ and $V_4$ over $X=M_0$.
Let us associate an $\N-$graded manifold $M$ to a graded vector bundle $\cV$, such that $\oM$ is equal to the total space of $\cV$, in all possible ways up to an isomoprhism. We will see that the latter is characterized by the 1st \v{C}ech cohomology group $H^1(M_0, S^2 V^*_2\otimes V_4)$. 

\noi By definition,
$\cO_{M_0}$ is always embedded into $\cO_M$ as functions of degree zero, while the vector bundle $V_2$ is uniquely determined by its dual $V_2^*$, whose sections coincide with functions on $M$ of degree $2$. Let $\xi_\a^a$ and $\psi_\a^\mu$ be local linear coordinates of degree $2$ and $4$ on $\cV$, defined on the coordinate chart $\uU_\a$, $(f_{\a\b}, h_{\a\b})$ be the transition cocycle for the vector bundle $V_2\oplus V_4$.
Any change of coordinates over the double overlaps $\uU_{\a\b}=\uU_\a\cap \uU_\b$ necessarily takes form
\beqn
\xi_\a^a=f_{\a\b}(x)^a_b \xi_\b^b\,, \hspace{3mm}
\psi_\a^\mu =h_{\a\b}(x)^\mu_\nu\psi_\b^\nu +
\frac{1}{2} \varphi_{\a\b}(x)^\mu_{ab}\psi^a_\a\psi^b_\a\, , 
\hspace{3mm} x\in \uU_{\a\b}\,,
\eeq
where $\varphi_{\a\b}(x)^\mu_{ab}$ are components of some $\varphi_{\a\b}\in \Gamma (\uU_{\a\b}, V_2^*\otimes V_4)$ in the basis, corresponding to $(\xi_\a^a, \psi_\a^\mu)$. The compatibility condition for the change of coordinates immediately implies
\beqn
\varphi_{\a\b} (x)=-\varphi_{\b\a}(x),\, \forall\, x\in \uU_{\a\b} \,\, and \,\, \varphi_{\a\ga}=\varphi_{\a\b}+\varphi_{\b\ga},\, \forall\, x\in \uU_{\a\b\ga}=\uU_\a\cap \uU_\b\cap \uU_\ga\,,
\eeq
which means that $\varphi$ is a 1-cocycle on $M_0$ with values in $S^2V_2^*\otimes V_4$.

\noi On the other hand, whenever we change degree $4$ coordinates from $\psi_\a^\mu$ to $\psi_\a^\mu+\frac{1}{2}\chi_{\a}(x)^\mu_{ab}\xi_\a^a\xi_\a^b$, $x\in \uU_\a$, the gluing cocycle changes from $\varphi_{\a\b}(x)$ to $\varphi_{\a\b}(x)+\chi_\b(x)-\chi_\a(x)$, $x\in \uU_{\a\b}$, i.e. the ambiguity is given by the cohomology class $[\varphi]\in H^1 (M_0, S^2 V^*_2\otimes V_4)$.

\noi In the smooth category $\varphi$ is always a coboundary, which means that there exists an isomorphism $M\simeq \oM$. In contrast to the smooth case, in the category of complex analytic manifolds it is no longer true.

\vskip 2mm\noindent
E.g. let $M_0=\mathbb{CP}^1$, $V_2=\cO(k)$ and $V_4 = \cO (l)$, where $\cO(q)$ denotes the holomorphic line bundle over the projective line with the Chern number $q$. It is unique up to an isomorphism of holomorphic line bundles and given by the $q-$th power of the dual to the canonical line bundle $\cL_{can}$ over $\mathbb{CP}^1$, which associates to a point of the projective line the corresponding 1-dimensional subspace of $\C^2$: 
$\cO(q) =\big(\cL_{can}^*\big)^{q}=\big(\cL_{can}\big)^{- q}$.

\noi Now $[\varphi]\in H^1 (\mathbb{CP}^1, \cO (l-2k))$. By Serre duality, the latter is isomorphic to $H^0 (\mathbb{CP}^1, \mathcal{K}\otimes S^2 V_2\otimes V^*_4)$, where $\mathcal{K}$ is the canonical class of the complex line.
Taking into account that, for $\mathbb{CP}^1$, $\mathcal{K}=\cO(-2)$,  we finally obtain
 \beqn
[\varphi]\in H^0 (\mathbb{CP}^1, \cO (2k-l-2))
\eeq
This cohomology group is non-trivial, i.e. the holomorphic bundle $\cO (2k-l-2)$ admits non-zero sections, if and only if $2k-l\ge 2$; as soon as it happens, by the above procedure the choice of such a holomorphic section determines a non-trivial twisting of $\oM$ into an $\N-$graded manifold $M$, which is not isomorphic to $\oM$.
\end{example}

\begin{example}[non-splitted $\Z-$manifold with splitted positive and negative parts]\label{ex:non-splittedZ} Assume that $M_+$ and $M_-$ are isomorphic to the total space of line bundles $\cV_{-1}$ and $\cV_{1}$ of degree $-1$ and $1$, respectively. To construct the whole $M$ out of this data, we should fix a gluing 1-cocycle $\varphi$ consisting of $\Z-$graded algebra isomorphisms $\varphi_{\a\b}$ over double overlaps $\uU_{\a\b}$. Let us choose a collection of local
nilpotent derivations of functions on the base with values in $\cV_{-1}^*\otimes\cV_{1}^*$, $\psi_{\a\b}\in \Gamma (\uU_{\a\b}, TM_0\otimes\cV_{-1}^*\otimes\cV_{1}^*)$ and put $\varphi_{\a\b} =\exp{\psi_{\a\b}}$. The corresponding change of coordinates takes the form $x^i \mto x^i + \left(\psi^i_{\a\b}\right)_{ab} (x) z^a w^b$, where $\psi_{\a\b} =\left(\psi^i_{\a\b}\right)_{ab} (x) z^a w^b\pt_{x^i}$ and
$x^i$, $z^a$ and $w^b$ are coordinates of degree $0$, $1$ and $-1$, respectively. 
\vskip 2mm\noindent
Now any $\psi$, the cohomology class of which $[\psi]\in H^1 (M_0, TM_0\otimes \cV^*_{-1}\otimes\cV^*_{1})$ is non-zero, gives us an obstruction to $M$ being isomorphic to the total space of $\cV=\cV_{1}\oplus \cV_{-1}$. While in the smooth case such an obstruction always vanishes (we used this in the proof of the Batchelor's theorem), in the complex analytic case it can be non-trivial.
\vskip 2mm\noindent
Similar to the Example \ref{ex:non-splittedN}, it is sufficient to take $M_0=\mathbb{CP}^1$, $\cV_1=\cO(k)[-1]$ and $\cV_{-1} = \cO (l)[1]$. Now 
\beqn 
 H^1 (\mathbb{CP}^1, T\mathbb{CP}^1\otimes\cO (-l-k)) \simeq H^0 (\mathbb{CP}^1, \mathcal{K}^2\otimes\cO (l+k))=
H^0 (\mathbb{CP}^1, \cO (l+k-4))
\eeq
since $\mathcal{K}=(T\mathbb{CP}^1)^*=\cO (-2)$. Whenever $l+k\ge 4$, a non-trivial twisting for $M$ exists.
\end{example}
\begin{rem}
 Note that in both examples the problems do not come from the infinite series issues that we discussed above: they are present already ``around'' the base manifold.  
\end{rem}

\section{Conclusion.}
In this paper we have started describing natural constructions for $\Z$-graded manifolds, namely 
defined the categorical properties of them.
As we have noticed, while the results were expected, the careful approach was definitely missing, it required a lot of clarifications and definitions. As mentioned in the introduction, the whole text is motivated by our previous work on integration of differential graded Lie algebras (\cite{DGLG}), where, again, the final result was natural, but the technicalities were highly non-trivial. With the construction of this paper one can now extend \cite{DGLG} to the honest $\Z$-graded case in a straightforward way by merely reformulating the Poincar\'e--Birkhoff--Witt theorem. We are also going to discuss some more advanced ``higher'' constructions on $\Z$-graded manifolds in \cite{cat-Z-fancy}.

\noi {\bf Acknowledgements.} \\
We thank Leonid Ryvkin and Camille Laurent-Gengoux for useful discussions. \\
 The meetings at early stages of this work were possible thanks to the grant no. 18-00496S of the Czech Science Foundation (for A.K.) and the career starting grant from the INSIS of CNRS (for V.S.). The work is also partially supported by the CNRS 80Prime project ``GraNum''.  
We are thankful to the Erwin Schr\"odinger International Institute for Mathematics and Physics (the program ``Geometry for Higher Spin Gravity: Conformal Structures, PDEs, and Q-manifolds'') as well as the Institut Henri Poincaré (Research in Paris program), that permitted the authors to gather in the same room and significantly advance on the manuscript.


 \appendix
\setcounter{equation}{0} 

\section{Direct limits and completions of algebras}
\label{app-filtr}
\noi Let $\sA$ and $\sC$ be two categories.

\begin{deff}[Direct and inverse limits]\label{def:limits_of_functor}
{\mbox{}\vskip 1mm}\noindent
An object $\cc\in\sC^0$ together with a collection of morphisms 
$F(\ca)\xrightarrow{\phi_\ca} \cc$ for all $\ca\in\sA^0$
is an inductive (direct) limit of a covariant functor 
$\sA\xrightarrow{F}\sC$ if:
 \begin{enumerate}
     \item For any morphism $\ca_1\xrightarrow{f}\ca_2$ in $\sA$ the following diagram is commutative:
     
     \beq\label{eq:diagram_direct_functor}
     \xymatrix{
      F(\ca_1) \ar[dr]_{\phi_{\ca_1}} \ar[rr]_{F(f)} && F (\ca_2)\ar[dl]^{\phi_{\ca_2}}\\
     & \cc &
             }
     \eeq
     \item For any other $\cc'\in\sC^0$ and $F(\ca)\xrightarrow{\phi'_\ca} \cc'$, which satisfy \eqref{eq:diagram_direct_functor}, there exists a unique morphism $\cc\to\cc'$, such that for any $\ca\in\sA^0$
     
     \beq\label{eq:universality_direct_limit}
     \xymatrix{
     & F(\ca) \ar[dr]^{\phi'_\ca}\ar[dl]_{\phi_\ca} & \\
     \cc \ar@{-->}[rr]^{\exists \,!} && \cc'
     }
     \eeq
 \end{enumerate}
 
 \noi An object $\cc\in\sC^0$ together with a collection of morphisms $\cc\xrightarrow{\psi_\ca} F(\ca)$ is a projective (inverse) limit of a contravariant functor $\sA\xrightarrow{F}\sC$ if the properties (1) and (2) above are satisfied with all arrows being reversed.
 
 \noi The standard notation for the direct (inverse) limit is $\cc=\varinjlim F$ and $\cc=\varprojlim F$, respectively.
\end{deff}

\begin{example}
$\sI$ is the category whose objects are natural numbers and arrows are pairs $(i\le j)$ with $i$ and $j$ being the source and the target, respectively. Then 
  \begin{itemize}
      \item a covariant functor $\sI\to\sC$ is a sequence of objects $\cc_i$, $i\in\N$ and morphisms $\cc_i\xrightarrow{f_{ij}}\cc_j$ for all $i\le j$, such that $f_{ik}=f_{jk}f_{ij}$ for all triples $i\le j\le k$. The definition of a direct limit $\big(\cc, \cc\xrightarrow{\phi_i}\cc_i\big)$ of $\big(\cc_i, f_{ij}\big)$, including the universality property, is encoded in the following commutative diagram:
      
      \beqn
      \xymatrix{
      \ca_i \ar[rr]^{f_{ij}} \ar[rd]^{\phi_i} \ar[rdd]_{\phi_i'} && \ca_j \ar[dl]_{\phi_j}\ar[ldd]^{\phi_j'}\\
     & \cc \ar@{-->}[d]^{} & \\
     & \cc' &
             }
      \eeq
      
      \item a contravariant functor $\sI\to\sC$ is a sequence of objects $\cc_i$, $i\in\N$ and morphisms $\cc_j\xrightarrow{p_{ij}}\cc_i$ for all $i\le j$, such that $p_{ik}=p_{ij}p_{jk}$ for all triples $i\le j\le k$. The inverse limit $\big(\cc, \cc_i\xrightarrow{\pi_i}\cc\big)$ of $\big(\cc_i, f_{ij}\big)$ will make the following diagram ommutative:
      
      \beqn
       \xymatrix{
       & \cc' \ar@{-->}[d] \ar[ldd]_{\pi_j'} \ar[rdd]^{\pi_i'} & \\
       & \cc \ar[ld]^{\pi_j} \ar[rd]_{\pi_i} & \\
       \ca_j \ar[rr]_{p_{ij}} && \ca_i
                 }
      \eeq
  \end{itemize}
\end{example}

\begin{example}\label{ex:lim_ass_alg}
$\sC$ is the category of associative rings, $\sA=\sI$ as in the previous example.
\begin{itemize}
    \item For a direct family $\big(R_i, f_{ij}\big)$ the explicit formula for the inductive limit of rings is
    $$\varinjlim \big(R_i, f_{ij}\big)=\Big(\coprod_{i\in\N} R_i\Big)\Big/\sim \,,$$
    where $\sim$ is the minimal equivalence relation generated by $a_i\sim f_{ij} \big(a_i\big)\in R_j$ for all $a_i\in R_i$. 
    \item Likewise, for an inverse family $\big(R_i, p_{ij}\big)$
    $$\varprojlim \big(R_i, p_{ij}\big)=
    \{ (a_1, a_2, \ldots )\in \prod_{i\in\Z} R_i \, |\, 
    a_i =p_{ij}\big(a_j\big)\,, \forall\, i\le j
    \}$$ 
\end{itemize}
\end{example}

\noi Let $R$ be an associative ring. By a decreasing filtration of $R$ we mean a sequence of ideals $F^p R$, $p\in\N$, such that
\beqn
R=F^0 R\supset F^1 R \supset F^p R\ldots \supset \{0\}
\eeq
In addition, we will require that 
\beq\label{eq:zero_intersection} \bigcap_{p\in\N}F^p R=\{0\}
\eeq

\noi Let us define a sequence of rings $R^{(p)}=R/F^p R$ for $p\in \N$. There are  natural surjective morphism of rings $p_{ij}\colon R^{(j)}\to R^{(i)}$ for all $i\le j$, which satisfies the property of an inverse (projective) family $p_{ik}=p_{jk}p_{ij}$ for all $i\le j\le k$. 

\noi Let $\widetilde{R}=\varprojlim R^{(p)}$; as it was mentioned in Example \ref{ex:lim_ass_alg}, $\widetilde{R}$ is identified with the sequence of elements $a^k\in R^{(k)}$, such that $p_{kl}(a^l)=a^k$ for all $k\le l$. The exists a canonical morphism of rings $\tilde{\phi}\colon R\to \widetilde{R}$, which associates to each $a\in R$ the corresponding equivalence class $a^p = a\!\mod F^p R$ for each $p$. It is clear that $(a^p)_{p\in N}$ satisfies the compatibility conditions written above. 

\noi By formula \eqref{eq:zero_intersection},  the morphism $\tilde{\phi}$ is injective. If, in addition, it is injective, we call $R$ complete with respect to the filtration. Otherwise $\widetilde{R}$ is called the formal completion of $R$.

\noi Define
\beqn
F^p \widetilde{R} =\{
\big( \underbrace{0, \ldots, 0}_p, a^{p+1}, \ldots \big)\in \widetilde{R}\, |\, p_{ij}(a^j)=0\,,\forall j>p, i\le p\, 
\} 
\eeq
It is a decreasing filtration of $\widetilde{R}$. Moreover, $\tilde{\phi}$ is a morphism of filtered rings, i.e. $\tilde{\phi}\big(F^p R\big)\subset F^p \widetilde{R}$ for all $p\in \N$. Obviously, $\widetilde{R}$ is complete with respect to this filtration.

\begin{example}
Let $R$ be the polynomial ring $K [x]$ for an associative ring $K$. There is a decreasing filtration by ideals $F^p R=x^p R$. The formal completion of $R$ is the ring of formal power series with coefficients in $K$: $\widetilde{K}=K [[x]]$. The morphism $\tilde{\phi}$ is the inclusion of polynomials into the ring of formal power series. It is injective, but not surjective.
\end{example}

 
 \section{Sheaves over groupoids}\label{sec:sheaf_over_groupoid} 
 
 \noi 
 In general, let $\cG \rightrightarrows \cG^0$ be a Lie groupoid with the source and target maps $s$ and $t$, respectively. By a sheaf on $\cG$ we mean a sheaf $\cS$ (of sets, groups, algebras or other algebraic structures) over 
$\cG^0$ together with an isomorphism of sheaves $\psi\colon s^*\cS\simeq t^*\cS$ over the ``space of arrows'' $\cG^1$, which is compatible with the multiplication $m\colon \cG^1\times_{\cG^0}\cG^1\to \cG^1$ in the next sense. Let $p_1$ and $p_2$ be the projections of the fibered product $\cG^1\times_{\cG^0}\cG^1$ onto the right and left factors, respectively. Then 
\beqn
m^* \psi = p_2^* (\psi) p_1^* \psi \colon (s\circ p_1)^*\cS \simeq 
(t\circ p_2)^*\cS
\eeq
This means that for any composible elements $g_1, g_2\in\cG^1$, such that $s (g_2)=t(g_1)$, one has 
\beqn
\psi_{g_2g_1}=\psi_{g_2}\psi_{g_1}\colon Stalk_{s(g_1)}\cS\simeq Stalk_{t(g_2)}\cS\,,
\eeq
where $Stalk_x\cS$ is the stalk of the sheaf $\cS$ at $x\in \cS^0$ ($\psi$ is viewed is a continuous section of the bundle of isomorphisms $s^*\big( Stalk\cS\big)\simeq t^* \big(Stalk\cS\big)$ in the discrete fiber topology).

\noi The following statements are true:
\begin{itemize}
    \item Assume that the quotient space $X=\cG^0/\cG$ is a smooth manifold. Then the $\cG-$invariant sections of $\cS$ 
    \beqn
    \cS^\cG =\{\cs\in \cS \, |\, t^* (\cs)=\psi \big(s^*(\cs)\big)\}
    \eeq
    gives rise to a sheaf on $X$.
    \item For any Morita equivalent groupoid $\cH$, there is a canonical sheaf on $\cH$ which produces the same sheaf on $X$ (under the above assumption). Therefore $\cS$ can be viewed as a sheaf on the corresponding stack.
\end{itemize}
To argue the last statement, let $P$ be a principal $\cG-\cH$-bibundle with the momentum maps $\mu_1$ and $\mu_2$ 
$$
\xymatrix{ & P\ar[dl]_{\mu_2}\ar[dr]^{\mu_1} & \\
\cG^0 && \cH^0
 }
$$
and $\cG-$action $\a\colon \cG^1\times_{\cG^0}P\to P$.
Then the sheaf over $\cH$ is given by $\cG-$invariant sections of the pull-back sheaf $\mu_2^*\cS$ over $P$:
$$
\big(\mu_2^*\cS\big)^{\cG} = \{\cs\in \mu_2^* \cS \,|\, 
\a^* (\cs)= p_2^* (\psi)p_1^* (\cs)\},\,
$$
where
$$
\xymatrix{ & \cG^1\times_{\cG^0}P\ar[dl]_{p_2}\ar[dr]^{p_1}  \\
\cG^1 & &  P
 }
$$
In particular, for a give open covering $\mathfrak{U}=(U_i)_{i\in I}$ of a topological space $X$, we define the \v{C}ech groupoid 
$\check{C}_{\mathfrak{U}}$, such that 
\beqn
\check{C}_{\mathfrak{U}}^0 = \coprod_{i\in I} U_i\, , \hspace{3mm}
\check{C}_{\mathfrak{U}}^1 = \coprod_{(i,j)\in I^2} U_{ij}
\eeq
where the source and the target maps are determined by the natural inclusions $U_{ij}\to U_j$ and $U_{ij}\to U_i$, respectively, the composition is given by the inclusion $U_{ij}\times_X U_{jk}\to U_{ik}$ and the inverse by the bijection $U_{ij}\to U_{ji}$. By construction of the \v{C}ech groupoid, it represents $X$, i.e.
\beqn
X= \check{C}_{\mathfrak{U}}^0 \, / \,\check{C}_{\mathfrak{U}}^1
\eeq
and the gluing cocycle for sheaves on the local sets $U_i$ can be interpreted as a sheaf on $\check{C}_{\mathfrak{U}}$. Any other open covering will produce a Morita equivalent groupoid (in this case the principal bibundle is the disjoint union of the open sets corresponding to the intersection of coverings). Therefore the result of the gluing procedure does not depend on the choice of a local covering.

\bibliography{BibGraded}

\end{document}